\documentclass[12 pt]{amsart}
\usepackage{amsmath,times,epsfig,amssymb,amsbsy,amscd,amsfonts,amstext,graphicx}
\usepackage{amscd,amssymb,amsmath}
\usepackage[arrow,matrix]{xy}
\usepackage{graphicx}

\theoremstyle{plain}
\newtheorem{thm}[subsection]{Theorem}
\newtheorem{lem}[subsection]{Lemma}
\newtheorem{prop}[subsection]{Proposition}
\newtheorem{cor}[subsection]{Corollary}

\theoremstyle{definition}
\newtheorem{rk}[subsection]{Remark}
\newtheorem{defn}[subsection]{Definition}
\newtheorem{ex}[subsection]{Example}

\numberwithin{equation}{section}
\newcommand{\M}{{\mathcal M}}

\newcommand{\A}{{\mathcal A}}

\newcommand{\CC}{{\mathcal C}}
\newcommand{\LL}{{\mathcal L}}

\newcommand{\Q}{\mathbb{Q}}
\newcommand{\R}{\mathcal{R}}
\newcommand{\C}{\mathbb{C}}

\newcommand{\FF}{\mathbb{F}}
\newcommand{\PP}{\mathbb{P}}
\newcommand{\RR}{\mathbb{R}}

\newcommand{\calQ}{{\mathcal Q}}
\newcommand{\FS}{{\mathcal{FS}}}
\DeclareMathOperator{\mult}{mult}

\begin{document}

\title [ connectivity of realization spaces of arrangements ]
{  On the connectivity of the  realization spaces of line arrangements }

\author{Shaheen Nazir}
\address{Department of Humanities and Sciences, 
National University of Computer \& Emerging Sciences, 
Lahore, PAKISTAN}
\email{shaheen.nazir@nu.edu.pk}


\author{Masahiko Yoshinaga}
\address{Department of Mathematics, 
Kyoto University, Kyoto, 606-8502, JAPAN}
\email{mhyo@math.kyoto-u.ac.jp}

\thanks{The second author has been supported by 
JSPS Grant-in-Aid for Young Scientists (B) 20740038. }

\dedicatory{}



\begin{abstract}
We prove that under certain combinatorial conditions, 
the realization spaces of line arrangements on the 
complex projective plane are connected. 
We also give several examples of arrangements with 
eight, nine and ten lines which have disconnected realization 
spaces. 
\end{abstract}

\maketitle

\section{Introduction}\label{intro}

Let $\A=\{H_1,\dots ,H_n\}$ be a line arrangement 
in the complex projective plane $\PP^2$ and denote by 
$M=M(\A)$, the corresponding arrangement complement. 
An arrangement $\A$ determines the incidence data $I(\A)$ 
(equivalently the intersection lattice $L(\A)$). 
This combinatorial data possesses the topological 
information, e.g. the cohomology algebra of $M$ are 
determined by the intersection lattice $L(\A)$ of $\A$. 
However, not all geometric information is determined by 
the incidence $I(\A)$. 
In 1993, Rybnikov \cite{Ry} posed an example of 
arrangements $\A_1$, $\A_2$ which have the same 
incidence but their fundamental groups 
are not isomorphic (see also \cite{ACCM1}). 
Nevertheless, in many cases the topological structures 
are determined by the combinatorial ones. They includes:
 \begin{enumerate}
\item 
Combining results of Fan \cite{Fa1},\cite{Fa2}, 
Garber, Teicher and Vishne \cite{GTV} and an unpublished 
work by Falk and Sturmfels (see \cite{cs}), 
if $n\leq 8$, then the fundamental group $\pi_1(M(\A))$ is 
determined by the combinatorics. 
\item 
In 2009, Nazir-Raza \cite{NR} introduced a complexity 
hierarchy of lattice: class $\CC_k$, and proved that 
if $\A$ is in $\CC_{\leq2}$, then the cohomology  
$H^*(M,\LL)$ with coefficients in a rank one local system 
$\LL$, is combinatorially determined.
\end{enumerate}
In this paper, we generalize these results by using 
the connectivity of the realization space $\R(I)$ of an 
incidence relation $I$. 
Indeed, the connectivity of realization spaces is 
related to the topology of the complements 
by Randell's lattice isotopy theorem. 

\begin{thm}
\label{thm:ra}
(Randell \cite{Ra}) 
If two arrangements are connected by a 
one-parameter family of arrangements 
which have the same lattice, then the 
complements are diffeomorphic, hence of 
the same homotopy type.
\end{thm}
Once the connectivity of the realization space 
$\R(I)$ is proved, then for any arrangements 
$\A_1$, $\A_2$ having the same incidences 
$I(\A_1)=I(\A_2)=I$, we can conclude that 
$M(\A_1)\cong M(\A_2)$ by Theorem \ref{thm:ra}. 
Since the realization space $\R(I)$ is a (quasi-projective) 
algebraic variety over $\C$, the irreducibility of 
$\R(I)$ implies the connectivity. 
(Note an irreducible algebraic variety is connected in 
the classical topology. For the proof, see 
\cite{Sh} chapter VII.) For our purposes, the 
following is useful. 
\begin{cor}
If $\R(I)$ is irreducible (in Zariski topology) and 
$I(\A_1)=I(\A_2)=I$, then $M(\A_1)\cong M(\A_2)$. 
\end{cor}

As far as the authors know, a systematic 
study of the connectivity of the realization space 
$\R(I)$ of line arrangements 
was initiated by Jiang and Yau \cite{JY} and 
subsequently by Wang and Yau \cite{WY}. 
They introduce the notion of graph associated 
to a line arrangement and under certain combinatorial 
conditions (``nice'' and ``simple'' arrangements), it is proved that 
$\R(I)$ is connected. In particular, the structure of 
fundamental groups are combinatorially determined. 
Explicit presentations for a class of combinatorially determined 
fundamental groups are also studied in \cite{EGT}. 

The purpose of this paper is to develop these ideas further. 
We will prove the connectivity of $\R(I)$ for ``inductively 
connected arrangement'' (Definition \ref{def:ic}) and 
``$\CC_{\leq 3}$ of simple type'' (Definition \ref{def:c3}). 
The relations between ``nice''(\cite{JY}), 
``simple''(\cite{WY}) and our classes are not clear at 
the moment. However up to $8$ lines, we will prove that 
all arrangements 
except for MacLane arrangement are contained in our class 
(\S\ref{sec:fg}, Proposition \ref{prop:8lines}). 
We also give a complete classification of disconnected 
realization space up to $9$ lines in \S\ref{sec:ex}.

\medskip

{\bf Acknowledgements.} 
The authors would like to thank 
the organizers of intensive research period 
``Configuration Spaces: Geometry, Combinatorics and Topology'' 
May--June 2010, at Centro di Ricerca Matematica Ennio De Giorgi, Pisa. 
The main part of this work was done during 
both authors were in Pisa. 
Our visits were supported by 
Centro di Giorgi, ICTP and JSPS. 
We also thank David Garber, Mario Salvetti 
and Filippo Callegaro for their useful 
comments to the previous version of this paper.

\section{Generality on the realization spaces of arrangements}
\label{sec:gen}
From now, we assume that $\A$ contains 
$H_i, H_j, H_k\in \A$ such that 
$H_i\cap H_j\cap H_k=\emptyset$ 
(thus excluding $n<3$ and pencil cases). 
Let $H_i\in \A$. $H_i$ is defined by 
$$
H_i=\{(x:y:z)\in \PP^2 \mid a_ix+b_iy+c_iz=0\}. 
$$
We may consider $(a_i:b_i:c_i)\in (\PP^2)^*$ 
as an element of dual projective plane. 
We call a triple $(H_i,H_j,H_k)$ an {\em intersecting 
triple} if $H_i\cap H_j\cap H_k\neq\emptyset$, or 
equivalently,
$$
\det(H_i,H_j,H_k):=\det\left(
  \begin{array}{ccc}
  a_i & b_i & c_i \\
   a_j & b_j & c_j \\
    a_k & b_k & c_k \\
    \end{array}
     \right)=0.
$$

\begin{defn}
Define the \textit{Incidence} of 
$\A$ by 
$$
I(\A):=\left\{\left.\{i,j,k\}\in 
\begin{pmatrix}
[n] \\
3 
\end{pmatrix}
\right| 
H_i\cap H_j\cap H_k\neq\emptyset\right\},
$$
where $ 
\begin{pmatrix}
[n]\\
3
\end{pmatrix}
=\left\{\{i,j,k\}\mid i,j,k\in \{1,2,\dots ,n\} 
\textrm{ mutually distinct}\right\}$.

\end{defn}
The set of all arrangements which have 
prescribed incidence $I$ is called 
the realization space of the incidence $I$. 
Let us define
$$
\R(I):=
\left\{
(H_1,\dots ,H_n)\in ((\PP^2)^*)^n \left| 
\begin{array}{l}
H_i\neq H_j\ \textrm{ for }i\neq j, \mbox{ and } \\
\det(H_i,H_j,H_k)=0 \mbox{ for } \{i, j, k\}\in I,\\
\det(H_i,H_j,H_k)\neq 0 \mbox{ for } \{i, j, k\}\notin I
\end{array}\right.
\right\}.
$$

It can be seen that $(H_1,\dots ,H_n)$ and $(gH_1,\dots ,gH_n)$ 
for $g\in PGL_3(\C)$ have the same incidence. Hence 
$PGL_3(\C)$ acts on $\R(I)$. Now, we will discuss 
the irreducibility of $\R(I)$.
\begin{defn}
Define 
$$
\overline{\R}(I):=
\left\{
(H_1,\dots ,H_n)\in ((\PP^2)^*)^n \left| 
\begin{array}{l}
H_i\neq H_j\ \textrm{ for }i\neq j, \mbox{ and } \\
\det(H_i,H_j,H_k)=0 \mbox{ for } \{i, j, k\}\in I
\end{array}\right.
\right\}.
$$
\end{defn}

\begin{ex}
Consider the incidence $I=\{\{1,2,3\}\}$ of 
$4$ lines $\{H_1, H_2, H_3, H_4\}$. 
Then 
$$
\R(I):=
\left\{
(H_1,\dots ,H_4)\in ((\PP^2)^*)^4 \left| 
\begin{array}{l}
H_i\neq H_j\ \textrm{ for }i\neq j, \mbox{ and } \\
\det(H_1, H_2, H_3)=0 \\
\det(H_1, H_2, H_4)\neq 0\\
\det(H_1, H_3, H_4)\neq 0\\
\det(H_2, H_3, H_4)\neq 0
\end{array}\right.
\right\}, 
$$
and, 
$$
\overline{\R}(I):=
\left\{
(H_1,\dots ,H_4)\in ((\PP^2)^*)^4 \left| 
\begin{array}{l}
H_i\neq H_j\ \textrm{ for }i\neq j, \mbox{ and } \\
\det(H_1, H_2, H_3)=0 
\end{array}\right. 
\right\}.
$$
\end{ex}

By definition, $\R(I)$ is a Zariski open subset of 
$\overline{\R}(I)$. Hence, $\overline{\R}(I)$ is 
irreducible implies that $\R(I)$ is irreducible and 
hence $\R(I)$ is connected (unless $\R(I)$ is empty). 

\begin{prop}
Assume that $\overline{\R}(I)$ is irreducible. 
Then $I=I(\A_1)=I(\A_2)$ implies that $M(\A_1)\cong M(\A_2)$. 
\end{prop}

\proof
From the assumption, $\R(I)$ is irreducible, hence 
connected. The result follows from Theorem \ref{thm:ra}. 
\endproof

\section{Connectivity and field of realization}
\label{sec:con}

In this section we establish several conditions on 
the incidence $I$ for the realization space $\R(I)$ 
to be connected. We also discuss field of definition, 
since in the case of $\leq 9$ lines, it is related to 
the connectivity of $\R(I)$. 

\begin{defn}
Let $\A$ be a line arrangement on $\PP_\C^2$. 
Denote by 
$$
\mult(\A)=\{
p\in\PP^2\mid
\mbox{$p$ is contained in $\geq 3$ lines of $\A$}\}. 
$$
We call $p\in\mult(\A)$ a multiple point. 
\end{defn}

The next lemma will be used frequently.
\begin{lem}\label{lem 1}
Let $\A=\{H_1,\dots ,H_n\}$ be a line arrangement in $\PP_\C^2$. 
Assume that $|H_n\cap\mult(\A)|\leq 2$. 
Set $\A'=\{H_1,\dots ,H_{n-1}\}$, $I=I(\A)$ and $I'=I(\A')$. 
If $\R(I')$ is irreducible, then $\R(I)$ is also irreducible. 
\end{lem}
\proof
Let $\mu=|H_n\cap\mult(\A)|$. 
By assumption $\mu\in\{0,1,2\}$. 
We claim that $\R(I)$ is a Zariski open subset of 
$\PP_\C^{2-\mu}$-fibration over $\R(I')$. 
Consider the projection $\pi:\R(I)\rightarrow \R(I')$ defined as
$(H_1,\dots ,H_n)\mapsto (H_1,\dots ,H_{n-1}).$
Let $p\in H_n\cap\mult(\A)$. 
Then $p$ is a (possibly normal crossing) intersection point of 
$\A'=\A\setminus H_n$. 

\medskip

\textbf{Case 1}: $\mu=2$.
Let $p_1,p_2\in H_n$ be multiple points of $\A$. 
In this case, $H_n$ can be  uniquely determined by 
$\A'$ as $H_n$ is the line connecting $p_1$ and $p_2$. 
Hence $\pi$ is an inclusion $\R(I)\hookrightarrow \R(I')$. 
The defining conditions of $\R(I)$ concerning $H_n$ 
other than ``$p_1, p_2\in H_n$'' 
are of the form $\det(H_i, H_j, H_n)\neq 0$, that is 
Zariski open conditions. Thus, in this case, 
$\pi:\R(I)\rightarrow \R(I')$ 
is a Zariski open embedding.

\medskip

\textbf{Case 2}: $\mu=1$. In this case, 
$H_n\cap\mult(\A)=\{p\}$. 
Suppose $p\in H_1,\dots ,H_t$ and $p\notin H_{t+1},\dots ,H_{n-1}$. 
Then the realization space can be described as 
$$
\R(I)=\left\{
(H',H_n)\in \R(I')\times (\PP^2)^*\left|
\begin{array}{l}
H_i\neq H_n, \mbox{ for }1\leq i\leq n-1,\\
\det(H_i, H_j, H_n)=0\ \textrm{ for } 1\leq i<j\leq t,\\
\det(H_i, H_j, H_n)\neq 0\ \textrm{ for others}
\end{array}
\right.\right\}.
$$
Note that the Zariski closed condition in the second line 
($\det(H_i,H_j,H_n)=0$) indicates that $H_n$ goes through 
$p=H_1\cap\dots\cap H_t$, which is equivalent to say that 
$H_n$ is contained in the dual projective line 
$p^{\perp}(\simeq \PP^1)\subseteq (\PP^2)^*$. 
Hence, $\R(I)$ is a Zariski open subset of 
$\PP^1$-fibration over $\R(I')$. 




\medskip

\textbf{Case 3}: $\mu=0$. 
In this case $H_n$ is generic to $\A'$. 
Hence $\R(I)$ is a Zariski open subset of 
$\R(I')\times(\PP^2)^*$. 
\endproof

Lemma \ref{lem 1} allows us to prove the irreducibility 
of $\R(I)$ by the inductive arguments. 

\begin{prop}
\label{prop:ind}
Let $\A=\{H_1, \dots, H_n\}$ be lines on $\PP_\C^2$. Define 
the subarrangements $\A_t=\{H_1, \dots, H_t\}$ 
($t=1, \dots, n$). If 
$|H_t\cap\mult(\A_t)|\leq 2$ for all $t$, then 
$\R(I(\A))$ is irreducible. 
\end{prop}
\proof
Induction on $t$ using Lemma \ref{lem 1}. 
\endproof

\begin{defn}
\label{def:ic}
A line arrangement $\A$ is said to be {\em inductively 
connected} (``i.c.'' for brevity) if there exists an 
appropriate numbering $\A=\{H_1, \dots, H_n\}$ of $\A$ 
which satisfies 
the assumption of Proposition \ref{prop:ind}. 
\end{defn}
Inductive 
connectedness is a combinatorial property. We also say 
the incidence $I=I(\A)$ is i.c. 
By Proposition \ref{prop:ind}, $R(I)$ is irreducible 
for i.c. incidence $I$. 

\begin{cor}
If $|\mult(\A)\cap H|\leq 2$ for all $H\in\A$, then 
$\A$ is i.c., hence 
$R(I(\A))$ is irreducible. 
\end{cor}

\begin{cor}
If $\R(I(\A))$ is disconnected, then there exists 
subarrangement $\A'\subset\A$ such that 
$$
|\mult(\A')\cap H|\geq 3, 
$$
for all $H\in\A'$. 
\end{cor}
\proof
If not, $\A$ is i.c. for any ordering. 
\endproof

\begin{rk}
It is easily seen that if the characteristic 
of the field is $\neq 2$ and 
$|\A|\leq 7$, every line arrangement is an i.c. arrangement. 
Obviously the set of all $\FF_2$-lines on $\PP_{\FF_2}^2$ is 
not i.c. In the case of characteristic zero, 
MacLane arrangement (Example \ref{ex:mc}) is the smallest 
one which is not i.c. 
\end{rk}

\begin{figure}[htbp]
\begin{picture}(400,214)(20,0)
\thicklines
\put(90,10){\footnotesize $H_3$}
\put(100,20){\line(0,1){180}}

\put(7,118){\footnotesize $H_8$}
\put(20,120){\line(1,0){160}}
\qbezier(180,120)(190,120)(200,100)

\put(7,78){\footnotesize $H_7$}
\put(20,80){\line(1,0){160}}
\qbezier(180,80)(190,80)(200,100)

\put(27,183){\footnotesize $H_9$}
\put(40,180){\line(1,-1){160}}

\put(27,17){\footnotesize $H_4$}
\put(40,20){\line(1,1){160}}

\put(27,37){\footnotesize $H_5$}
\put(40,40){\line(1,1){130}}
\qbezier(170,170)(180,180)(200,180)

\put(110,10){\footnotesize $H_2$}
\put(120,20){\line(0,1){160}}
\qbezier(120,180)(120,190)(100,200)

\put(7,98){\footnotesize $H_6$}
\put(20,100){\line(1,0){180}}

\put(15,205){\footnotesize $H_1$}
\put(200,0){\line(0,1){190}}
\qbezier(200,190)(200,200)(190,200)
\put(20,200){\line(1,0){170}}


\put(300,20){\line(0,1){180}}

\put(220,80){\line(1,0){160}}
\qbezier(380,80)(390,80)(400,100)

\put(220,120){\line(1,0){160}}
\qbezier(380,120)(390,120)(400,100)

\put(240,180){\line(1,-1){160}}

\put(240,160){\line(1,-1){130}}
\qbezier(370,30)(380,20)(400,20)

\put(240,20){\line(1,1){160}}

\put(240,40){\line(1,1){130}}
\qbezier(370,170)(380,180)(400,180)

\put(320,20){\line(0,1){160}}
\qbezier(320,180)(320,190)(300,200)

\put(220,100){\line(1,0){180}}

\put(400,0){\line(0,1){190}}
\qbezier(400,190)(400,200)(390,200)
\put(220,200){\line(1,0){170}}

\end{picture}
      \caption{An i.c. arrangement $\A_1$ and non i.c. 
arrangement $\A_2$. Both are $\CC_3$ of simple type.}
\label{fig:ic}
\end{figure}
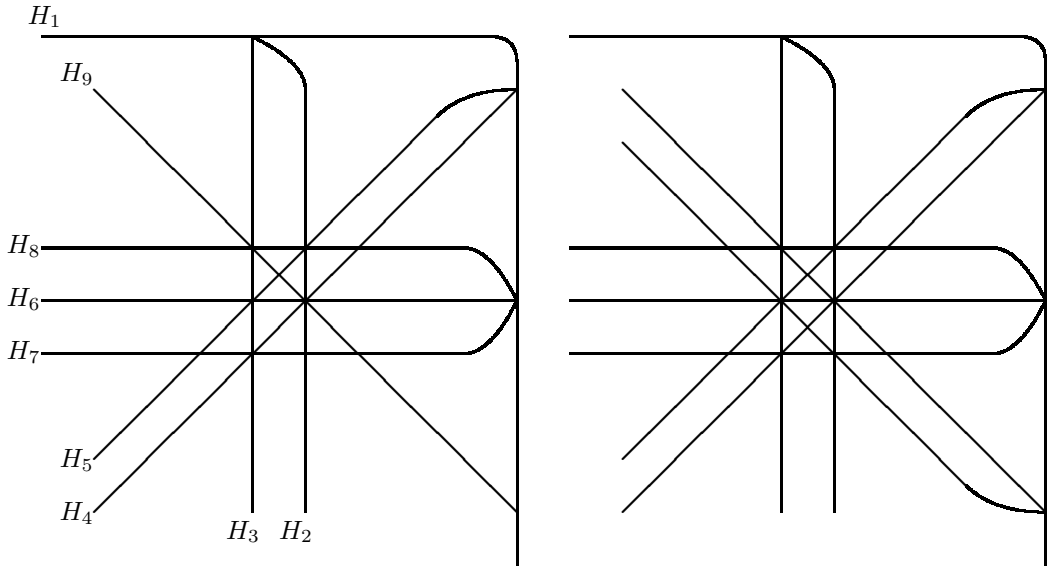

\begin{ex}
Let $\A_1$ (resp. $\A_2$) be a line arrangement 
defined as left of Figure \ref{fig:ic} (resp. right). 
Then $\A_1$ is i.c., but $\A_2$ is not i.c. 
(Each line $H\in\A_2$ has at least $3$ multiple 
points.) 
\end{ex}

Let $K\subset \C$ be a subfield, $I$ an incidence. 
The incidence $I$ is realizable over the field $K$ 
if the the set of $K$-valued points $\R(I)(K)$ is nonempty. 
(Equivalently, there exists an arrangement $\A$ with the 
coefficients of defining linear forms in 
$K$ satisfying $I=I(\A)$.) 
The next Lemma can also be proved similarly as Lemma \ref{lem 1}. 

\begin{prop}
\label{prop2}
With notation as in Lemma \ref{lem 1}, 
if the $K$-valued points $\R(I')(K)$ is 
Zariski dense in $\R(I')(\C)$, then $\R(I)(K)$ 
is Zariski dense in $\R(I)(\C)$. 
In particular, $\R(I)(K)\neq \emptyset$ and 
$I$ is  realizable over $K$. Every i.c. arrangement 
is realizable over $\Q$. 
\end{prop}

Next we discuss connectivity of $\R(I)$ for another 
type of incidence. 
\begin{defn}
Let $k$ be a non-negative integer. We say that a 
line arrangement $\A$ (or its incidence  $I(\A)$) 
is of type $\CC_k$ if $k$ is the minimal number of 
lines in $\A$ 
containing all the multiple points. 
\end{defn}
For instance $k=0$ corresponds to nodal arrangements, 
while $k=1$ corresponds to the case of a nodal 
{\it affine} arrangement. 
Note that $k=k(\A)$ is combinatorially defined, i.e. 
depends only on the intersection lattice $L(\A)$. 

\begin{thm}\label{th 1}
Let $\A=\{H_1,\dots ,H_n\}$ be a line arrangement in 
$\PP^2_{\C}$ of class $\CC_{\leq 2}$ (i.e., either $\CC_0, \CC_1$ 
or $\CC_2$). 
Then $\A$ is i.c. In particular, 
the realization space $\R(I(\A))$ is irreducible. 
\end{thm}

\proof
By assumption, we may say that all multiple 
points are on $H_1\cup H_2$. 
For $i\geq 3$, as $|H_i\cap (H_1\cup H_2)|\leq 2$, 
there are at most two multiple points on $H_i$. 
Hence the subarrangements 
$\A_t:=\{H_1, \dots, H_t\}$ 
($t=1,\dots, n$) 
satisfy the assumption 
of Proposition \ref{prop:ind}. Thus $\R(I(\A))$ 
is irreducible. 
\endproof

\begin{rk}
Under the assumption of Theorem \ref{th 1}, 
using Proposition \ref{prop2}, 
we can prove that $I(\A)$ is realizable over $\Q$. 
\end{rk}
The irreducibility of the realization spaces are not 
guaranteed for class $\CC_3$ in general 
(see Example \ref{ex:mc+}). 
Now we introduce a subclass of $\CC_3$.
\begin{defn}
\label{def:c3}
Let $\A$ be an arrangement of type $\CC_3$.
Then $\A$ is called $\CC_3$ of simple type 
if there are $H_1, H_2, H_3\in \A$ such that all 
multiple points are in $H_1\cup H_2 \cup H_3$ and 
one of the following holds: 
\begin{description}
\item[(i)] 
$H_1\cap H_2 \cap H_3=\emptyset$ and there is only one 
multiple point on $H_1\setminus (H_2\cup H_3)$;
\item[(ii)] $H_1\cap H_2 \cap H_3\neq \emptyset$. 
\end{description}

\begin{figure}[htbp]
  \begin{picture}(360,150)(55,35)
 \put(70,60){\line(1,1){100}}
 \put(230,60){\line(-1,1){100}}
 \put(70,80){\line(1,0){160}}
  \put(75,85){$H_2$}
  \put(125, 78){$\bullet$}
  \put(150,78){$\bullet$}
\put(175,110){$\bullet$}
 \put(110,160){$H_{1}$}
 \put(175,160){$H_{3}$}
  \put(135, 125){$\bullet$}
  \put(120,110){$\bullet$}
\put(109,100){$\bullet$}
\put(137,40){(i)}

  \put(270,60){\line(1,1){100}}
 \put(430,60){\line(-1,1){100}}
  \put(350,60){\line(0,1){100}}

 \put(370,163){$H_{3}$}
 \put(325,163){$H_2$}

\put(347,138){$\bullet$}
 \put(348,163){$H_{1}$}
 \put(347, 95){$\bullet$}
   \put(347, 75){$\bullet$}
  \put(347,110){$\bullet$}
\put(309,100){$\bullet$}
 \put(320, 111){$\bullet$}
  \put(373,110){$\bullet$}
\put(393,90){$\bullet$}
 \put(342,40){(ii)}


\end{picture}
\caption{$\CC_3$ of simple type}
\end{figure}
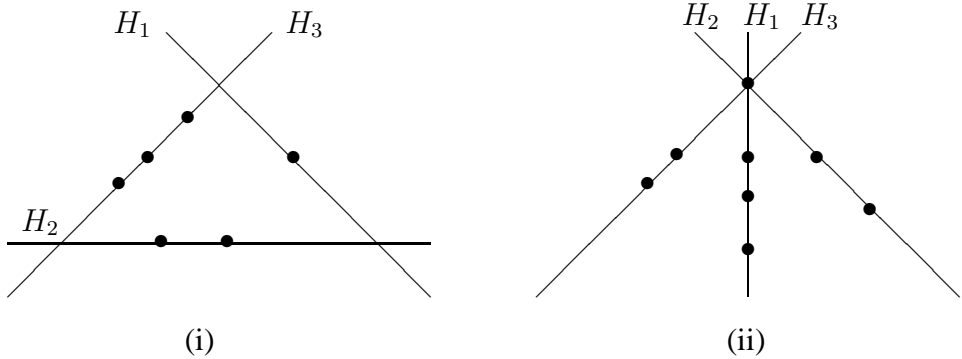
\end{defn}

\begin{ex}
The both line arrangements defined in 
Figure \ref{fig:ic} are $\CC_3$ of simple type. 
(E.g. $\mult(\A)\subset H_1\cup H_2\cup H_3$.) 
\end{ex}

\begin{thm}\label{th 2}
Let $\A$ be  an arrangement of 
$\CC_3$ of simple type. Then $\R(I(\A))$ is irreducible.

\end{thm}
\proof
The proof is divided into two parts according to (i) and (ii) of the definition of $\CC_3$ of simple type.

\subsection*{Case: (i)}
By the assumption, there exist $H_1, H_2, H_3\in \A$ 
which satisfy the condition (i). 
Let $p\in H_1\setminus (H_2\cup H_3)$ be the unique multiple point. 
Let us assume that $H_4,\dots ,H_t$ contain $p$ and 
$H_{t+1},\dots ,H_n$ do not contain $p$. 
For $i\geq {t+1} $, $H_i$ has at most two multiple points. 
By Lemma \ref{lem 1}, it suffices to prove the 
irreducibility for $\A'=\{H_1,\dots ,H_t\}$. 
However in this case, there are at most two multiple points: 
one is $p$ and the other possibility is $H_2\cap H_3$. 
Hence by Theorem \ref{th 1}, $\R(I(A'))$ is irreducible 
and so is $\R(I(\A))$.

\subsection*{Case: (ii)}

By the assumption, there exist $H_1, H_2, H_3\in \A$ 
which satisfy the condition (ii) of the definition. 
Let $O=H_1\cap H_2 \cap H_3$. If $H_i$  ($i\geq 4$) 
passes through $O$, then there is only one multiple 
point on $H_i$. Thus, by Lemma \ref{lem 1}, 
the irreducibility of $\R(I(\A))$ is reduced to $\R(I(\A'))$, 
where $\A'=H_1\cup H_2\cup H_3\cup\bigcup_{O\notin H_j}H_j$. 
We shall prove the irreducibility of $\R(I(\A'))$ by 
describing $\overline{\R}(I(\A))/PGL_3(\C)$ explicitly. 
By the $PGL_3(\C)$-action, we may fix as follows: 
$H_1=\{(x:y:z) \mid x=0\}, H_2=\{(x:y:z) \mid x=z\}$ and 
$H_3=\{(x:y:z) \mid z=0\}$, so $O=H_1\cap H_2\cap H_3=(0:1:0)$. 
We list all intersections on $H_i\setminus\{O\}$, (i=1, 2, 3): 
$$
P_\alpha(0:a_\alpha:1)\in H_1, (\alpha=1,\dots ,r, a_\alpha\in \C),
$$
$$
Q_\beta(1:b_\beta:1)\in H_2, (\beta=1,\dots ,s, b_\beta\in \C),
$$
$$
R_\gamma(1:c_\gamma:0)\in H_3, (\gamma=1,\dots ,t, c_\gamma\in \C).
$$
Every line $H_i$ ($i\geq 4$) in $\A'$, 
can be described as a line connecting 
$P_{\alpha_i}$ and $Q_{\beta_j}$. Hence, the quotient space
$\R(I(\A))/PGL_3(\C)$ 
can be embedded in the space 
$\C^{r+s+t}=\{(a_{\alpha},b_{\beta}, c_{\gamma})\}$. 
(More precisely, here we consider 
$X:=\C\times\C^*\times\R(I(\A))/PGL_3(\C)$. 
Because we fix only $H_1, H_2, H_3$ and the isotropy 
subgroup is $\{g\in PGL_3(\C)\mid gH_i=H_i, i=1, 2, 3\}
\simeq\C\times\C^*$.) 
Thus, we can describe the realization space by using 
the parameters $a_{\alpha},b_{\beta}, c_{\gamma}$.

Suppose $H_i$ $(i\geq 4)$ passes through 
$P_{\alpha_i}, Q_{\beta_i}, R_{\gamma_i}$. 
These three points are collinear if and only if 
$$
\det\left(
     \begin{array}{ccc}
       0 & a_{\alpha_i} &1 \\
      1 & b_{\beta_i} & 1 \\
       1 & c_{\gamma_i} & 0 \\
     \end{array}
   \right)=a_{\alpha_i}-b_{\beta_i}+c_{\gamma_i}=0.
$$
Collecting these linear equations together, we have
$$
(a_1, \dots, a_r, b_1, \dots, b_s, c_1, \dots, c_t)\cdot A=
\left(
\begin{array}{c}
0 \\
\vdots \\
0 \\
\end{array}
\right),
$$
where $A$ is a $(r+s+t)\times (n-3)$ matrix with entries 
$\pm1$ or $0$. Thus the space $X$ can be described as 
$$
X=\left\{
(a_{\alpha},b_{\beta}, c_{\gamma})\in \C^{r+s+t}\left| 
\begin{array}{l}
(a_{\alpha},b_{\beta}, c_{\gamma})\cdot A=0, \\
a_{\alpha}\neq a_{\alpha'},b_{\beta}\neq b_{\beta'}, 
c_{\gamma}\neq c_{\gamma'}, \\
\textrm{and other Zariski open conditions.} 
\end{array}
\right.\right\}. 
$$
Since $\ker A$ is isomorphic to $\C^K$ for some $K\geq 0$, 
the Zariski open subset $X\subset\C^K$ 
is irreducible.
\endproof

Thus we have proved that if $\A$ is either in the 
class $\CC_{\leq 2}$ or $\CC_3$ of simple type 
(``{\em $\CC_{\leq 3}$ of simple type}'' for short), 
$\R(I(\A))$ is connected. As it is mentioned, 
there are 
arrangements in $\CC_3$ of non-simple type which have 
disconnected realization spaces (Example \ref{ex:mc+}). 

By lattice isotopy theorem, we have 
\begin{cor}
Let $\A_1, \A_2$ be arrangements in $\PP^2$ 
of $\CC_{\leq 3}$ of simple type. If $I(\A_1)=I(\A_2)$, 
then the pairs $(\PP^2,\cup_{H\in \A_1}H)$ and 
$(\PP^2,\cup_{H\in \A_2}H)$ are homeomorphic. 
\end{cor}

\begin{rk}
Under the assumption of Theorem \ref{th 2}, 
$\R(I(\A))(\Q)$ is Zariski dense in $\R(I(\A))(\C)$, 
hence realizable over $\Q$. 
The proof is similar. Case (i) uses Proposition \ref{prop2} 
and in case (ii), we note that the matrix $A$ is with 
$\Q$-coefficients. Hence
$\ker A$ has $\C$-valued points if and only 
if it has $\Q$-valued points.

\end{rk}

\section{Application to the fundamental groups}
\label{sec:fg}

In this section, as an application of the 
connectivity theorem, we prove the following:
\begin{thm}\label{th 3}
Let $\A_1$ and $\A_2$ be two line arrangements in $\PP^2_{\C}$.  
Suppose that $|\A_1|=|\A_2|\leq 8$ and $I(\A_1)=I(\A_2)$. 
Then $$(\PP^2_{\C},\A_1)\cong (\PP^2_{\C},\A_2).$$
\end{thm}

\begin{cor}
Under the same assumption, we have 
$$\pi_1(M(\A_1))\simeq\pi_1(M(\A_2)).$$
\end{cor}

Thus the isomorphism classes of the fundamental 
groups are combinatorial for $n\leq 8$. 

The proof is done by using Theorem \ref{th 2} in \S\ref{sec:con}. 
Indeed, for almost all cases, $\A$ is of class $\CC_{\leq 3}$ 
of simple type. 
Hence the realization space is connected. 
However there is exception (unique up to the $PGL$-action and 
the complex conjugation).

\begin{ex} (MacLane arrangement $\M^\pm$) 
\label{ex:mc}
Let $\omega_{\pm}:=\frac{1\pm\sqrt{-3}}{2}$ 
be the roots of the quadratic equation $x^2-x+1=0$. 
Consider $8$ lines $\M^\pm=\{H_1,\dots ,H_8\}$ defined by: 
$$
\begin{array}{lll}
  H_1:x=0, & H_2:x=z, & H_3:x=\omega_{\pm}z, \\
  H_4:y=0, & H_5:y=z, & H_6:y=\omega_{\pm}z, \\
  H_7:x=y, & H_8:\omega_{\pm}x+y=\omega_{\pm}. &
\end{array}
$$

\begin{figure}[htbp]
\begin{picture}(200,170)(0,0)
\thicklines

\put(50,0){\line(0,1){140}}
\put(90,0){\line(0,1){165}}
\put(130,0){\line(0,1){140}}

\put(50,140){\line(2,1){50}}
\put(130,140){\line(-2,1){50}}

\put(30,10){\line(1,0){170}}
\put(30,50){\line(1,0){195}}
\put(30,90){\line(1,0){170}}

\put(200,10){\line(1,2){25}}
\put(200,90){\line(1,-2){25}}

\put(40,0){\line(1,1){170}}

\put(80,0){\line(1,1){60}}
\put(30,70){\line(1,1){30}}
\multiput(60,100)(5,2){21}{\circle*{1}}
\put(160,140){\line(1,-1){20}}
\multiput(140,60)(2,3){21}{\circle*{1}}

\put(46,-10){$H_1$}
\put(86,-10){$H_2$}
\put(126,-10){$H_3$}

\put(13,10){$H_4$}
\put(13,50){$H_5$}
\put(13,90){$H_6$}

\put(212,162){$H_7$}

\put(13,70){$H_8$}
\put(66,-10){$H_8$}


\end{picture}
      \caption{MacLane Arrangement $\M^\pm$}
\label{fig:maclane}
\end{figure}
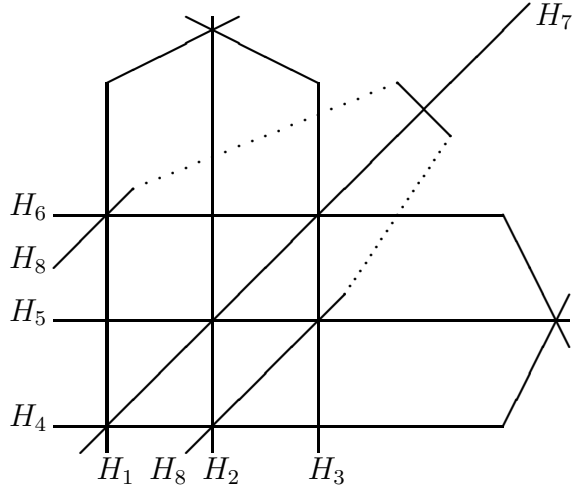
The MacLane arrangement is not of type $\CC_{\leq 3}$, 
but of type $\CC_4$ 
(e.g. all multiple points are contained in 
$H_1\cup H_2\cup H_3\cup H_4$), 
and the realization space has two connected components. 
$$\R(I)/PGL_3(\C)=\{\M^+,\M^-\}.$$
However the corresponding complements 
$M(\M^+)$ and $M(\M^-)$ are diffeomorphic 
under complex conjugation. 
Hence the complements have isomorphic fundamental groups.
\end{ex}

To prove Theorem \ref{th 3}, it is suffices to 
prove the following. 
\begin{enumerate}
  \item If $n\leq 5$, then $\A$ is in class $\CC_{\leq1}$;
  \item $n\leq 6$, then $\A$ is in class $\CC_{\leq 2}$;
  \item $n\leq 7$, then $\A$ is in class $\CC_{\leq 3}$ of simple type;
  \item $n=8$, then $\A$ is  either in class  $\CC_{\leq 3}$ 
of simple type or isomorphic to the MacLane arrangement $\M^{\pm}$.
\end{enumerate}

\textit{Proof of (1) and (2):}
\begin{enumerate}
  \item If a line arrangement  is in class $\CC_2$, 
then it is clear that there should be at least six lines. 
Thus, for $n\leq 5$, $\A$ is in class $\CC_1$. 

  \item Let $H\in \A$ and $\A'=\A\setminus\{H\}$. 
Then by $(1)$, there is a line $H'\in \A'$ such that 
all multiple points of $\A'$ are contained in $H'$, 
therefore, all multiple points of $\A$ are  in $H\cup H'$. 
Thus, $(2)$ holds.

\end{enumerate}

The following is the key lemma for our classification. 

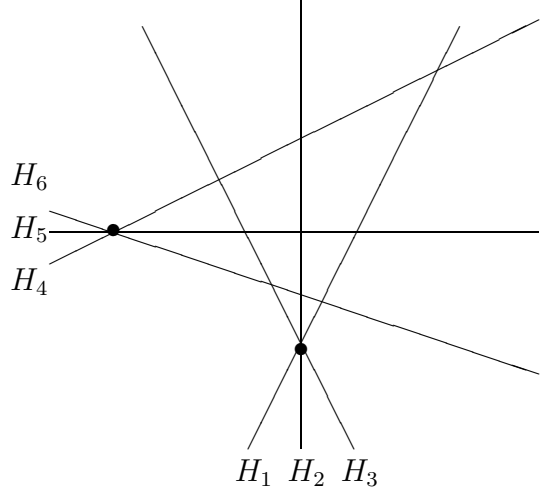
\begin{figure}[h]
\begin{center}
\begin{picture}(360,140)

  \put(180,0){\line(0,2){170}}
 \put(200,0){\line(-1,2){80}}
 \put(160,0){\line(1,2){80}}
 \put(85,70){\line(2,1){185}}
 \put(85,90){\line(3,-1){185}}
  \put(85,82){\line(1,0){185}}


 \put(155,-12){$H_1$}
 \put(175,-12){$H_2$}
 \put(195,-12){$H_3$}

 \put(70,60){$H_4$}
 \put(70,80){$H_5$}
 \put(70,100){$H_6$}

\put(106,80){$\bullet$}
\put(177, 35){$\bullet$}
\end{picture}
\caption{$6$ lines contained in a non $\CC_{\leq 3}$-simple type $\A$}
\label{fig:6lines}
\end{center}
\end{figure}

\begin{lem}\label{lem 3}
Let $\A$ be a line arrangement which is not in class $\CC_{\leq 3}$ 
of simple type. Then there exist $H_1, H_2,\dots ,H_6\in \A$ 
satisfying
$H_1\cap H_2 \cap H_3\neq \emptyset$, 
$H_4\cap H_5\cap H_6\neq \emptyset$, and 
$(H_1\cup H_2\cup H_3)\cap (H_4\cup H_5 \cup H_6)$ 
consists of $9$ points. (Figure \ref{fig:6lines}.)
\end{lem}
\proof
Suppose $H_1\cap H_2 \cap H_3=\{p\}\neq \emptyset$. 
Then there exists a multiple point which is not 
contained in $H_1\cup H_2 \cup H_3$, otherwise, 
$\A$ will be in class  $\CC_{\leq 3}$ of simple type. 
Suppose $H_4\cap H_5\cap H_6=\{q\}\neq \emptyset$ be such a 
multiple point. 
If there is no line which passes $p$ and $q$, then 
$H_1, \dots, H_6$ satisfy the conditions. 
If there exists a multiple point of multiplicity $3$, 
say $p$, then $H_4, H_5$ and $H_6$ do not pass $p$. 
Then again $H_1, \dots, H_6$ satisfy the conditions. 
If both $p$ and $q$ have multiplicity $\geq 4$ and 
there is a line passing $p$ and $q$, then there exist 
lines $H_1, \dots, H_7$ such that 
$\{p\}=H_1\cap H_2\cap H_3\cap H_4$ and 
$\{p\}=H_4\cap H_5\cap H_6\cap H_7$. Then 
$H_1, H_2, H_3, H_5, H_6, H_7$ satisfy 
the conditions. 
\endproof

\begin{prop}
Let $\A$ be a line arrangement with $|\A|=7$. 
Then $\A$ is in class $\CC_{\leq 3}$ of simple type.

\end{prop}

\proof
Suppose that $\A$ is not in class $\CC_{\leq 3}$  
of simple type. Then there exist $6$ lines $H_1,\dots ,H_6\in \A$ 
satisfying the conditions of Lemma \ref{lem 3}. 
So, all multiple points of $\A$ 
are either  $H_1\cap H_2 \cap H_3$, $H_4\cap H_5\cap H_6$ 
or contained in the line $H_7$.

%
%
%
%
%
%
%
%
Hence, all multiple points are contained in $H_1\cup H_4\cup H_7$. 
Moreover, as multiple points on $H_1\setminus (H_4\cup H_7)$ 
are at most one, $\A$ is in $\CC_{\leq 3} $ of simple type, which 
is a contradiction. 
\endproof

\begin{prop}
\label{prop:8lines}
Let $\A$ be a line arrangement with $|\A|=8$. 
Then $\A$ is either in class $\CC_{\leq 3}$ of simple type 
or $\A=\M^{\pm}$, the MacLane arrangement. 
\end{prop}
\proof
Suppose that $\A$ is not in class $\CC_{\leq 3}$ of simple type. 
Then by Lemma \ref{lem 3}, we have six lines 
$L_1, L_2, L_3, K_1, K_2, K_3\in \A$ such that 
\begin{itemize}
\item $L_1\cap L_2\cap L_3\neq \emptyset$, 
$K_1\cap K_2\cap K_3\neq \emptyset$, and 
\item Let $Q_{ij}:=L_i\cap K_j$. Then  
$Q_{ij}=Q_{i'j'}$ only if $i=i',j=j'$. 
\end{itemize}
Let us denote by $\calQ:=\{Q_{ij}\mid i,j=1,2,3\}$ 
the set of $9$ intersections of 
$(L_1\cup L_2\cup L_3)\cap(K_1\cup K_2\cup K_3)$. 
Suppose $\A=\{L_1, L_2, L_3, K_1, K_2, K_3, H_7, H_8\}$. 
We divide the cases according to the cardinality of 
$H_7\cap \calQ$ and $H_8\cap \calQ$. We may assume that 
$0\leq |H_7\cap \calQ|\leq |H_8\cap \calQ|\leq 3$. 

\medskip

\begin{figure}[htbp]
\begin{picture}(300,143)(0,0)
\thicklines

\put(50,50){\line(2,1){160}}
\put(20,10){\line(2,1){220}}
\put(100,25){\line(2,1){170}}

\put(100,25){\line(-1,0){70}}
\put(50,50){\line(0,-1){40}}

\put(213,133){$L_1$}
\put(243,123){$L_2$}
\put(273,113){$L_3$}

\put(250,50){\line(-2,1){160}}
\put(280,10){\line(-2,1){220}}
\put(200,25){\line(-2,1){170}}

\put(200,25){\line(1,0){70}}
\put(250,50){\line(0,-1){40}}

\put(15,113){$K_1$}
\put(45,123){$K_2$}
\put(75,133){$K_3$}

\put(100,75){\circle*{4}}
\put(94.5,82.5){\footnotesize $Q_{11}$}

\put(125,87.5){\circle*{4}}
\put(119.5,95){\footnotesize $Q_{12}$}

\put(150,100){\circle*{4}}
\put(144.5,107.5){\footnotesize $Q_{13}$}

\put(165,10){\line(0,1){130}}
\put(160,0){$H_7$}

\multiput(125,62.5)(25,12.5){3}{\circle*{4}}
\multiput(150,50)(25,12.5){3}{\circle*{4}}

\end{picture}
      \caption{Case $1$: $\calQ\cap H_7=\emptyset$.}
\label{fig:case1}
\end{figure}
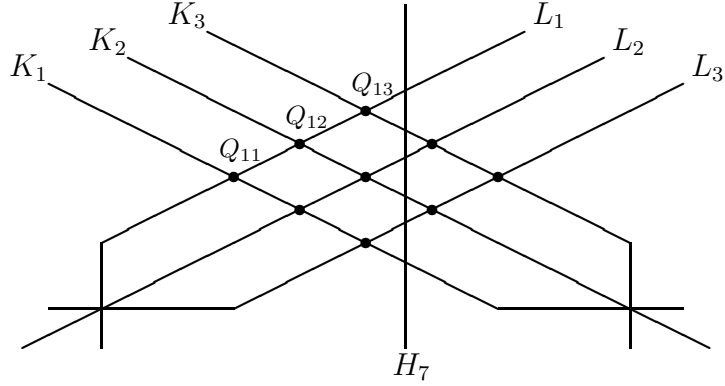

\noindent
{\bf Case 1}: $|H_7\cap \calQ|=0$ (Fig. \ref{fig:case1}). 
In this case, every multiple point of $\A$ is 
contained in $K_1\cup L_1\cup H_8$ and there 
are at most one multiple point in 
$K_1\setminus (L_1\cup H_8)$. 
Hence, $\A$ is in $\CC_{\leq 3}$ of simple type. 

%
%
%
%
%
%

\medskip

\begin{figure}[htbp]
\begin{picture}(300,143)(0,0)
\thicklines

\put(50,50){\line(2,1){160}}
\put(20,10){\line(2,1){220}}
\put(100,25){\line(2,1){170}}

\put(100,25){\line(-1,0){70}}
\put(50,50){\line(0,-1){40}}

\put(213,133){$L_1$}
\put(243,123){$L_2$}
\put(273,113){$L_3$}

\put(250,50){\line(-2,1){160}}
\put(280,10){\line(-2,1){220}}
\put(200,25){\line(-2,1){170}}

\put(200,25){\line(1,0){70}}
\put(250,50){\line(0,-1){40}}

\put(15,113){$K_1$}
\put(45,123){$K_2$}
\put(75,133){$K_3$}

\put(100,75){\circle*{4}}

\put(125,87.5){\circle*{4}}

\put(150,100){\circle*{4}}

\put(162.5,10){\line(1,4){30}}
\put(160,0){$H_7$}

\multiput(125,62.5)(25,12.5){3}{\circle*{4}}
\multiput(150,50)(25,12.5){3}{\circle*{4}}

\end{picture}
      \caption{Case $2$: $\calQ\cap H_7=\{Q_{32}\}$.}
\label{fig:case2}
\end{figure}
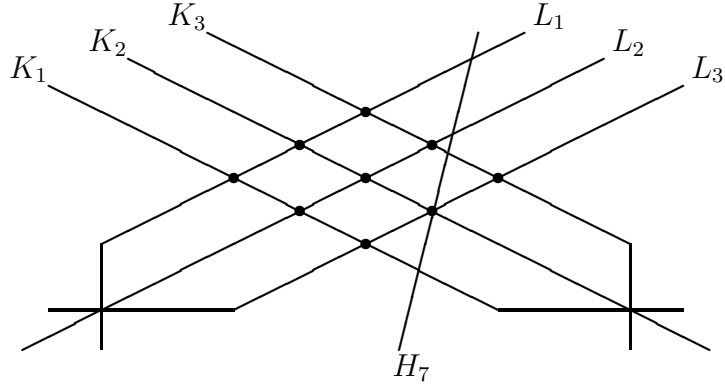

\noindent
{\bf Case 2}: $|H_7\cap Q|=1$ (Fig. \ref{fig:case2}). 
Let $H_7\cap \calQ =L_i\cap K_j=\{Q_{ij}\}$.  
Then every multiple point of $\A$ is contained in 
$K_j\cup L_i\cup H_8$ and there are at most 
one multiple point in $K_j\setminus (L_i\cup H_8)$. 
Hence, $\A$ is in $\CC_3$ of simple type. 

%
%
%
%

\medskip

The rest cases are $2\leq |H_7\cap \calQ|\leq |H_8\cap \calQ|\leq 3.$

\medskip

\begin{figure}[htbp]
\begin{picture}(300,150)(0,0)
\thicklines

\multiput(75,0)(75,0){3}{\line(0,1){150}}
\multiput(40,25)(0,50){3}{\line(1,0){220}}

\put(45,5){\line(3,2){210}}

\put(55,-10){\footnotesize $L_1:x=0$}
\put(130,-10){\footnotesize $L_2:x=1$}
\put(205,-10){\footnotesize $L_3:x=t$}

\put(-10,22){\footnotesize $K_1:y=0$}
\put(-10,72){\footnotesize $K_2:y=1$}
\put(-10,122){\footnotesize $K_3:y=t$}

\put(0,-2){\footnotesize $H_8:x=y$}

\put(78,14){\footnotesize $Q_{11}$}
\put(153,14){\footnotesize $Q_{21}$}
\put(228,14){\footnotesize $Q_{31}$}

\put(78,64){\footnotesize $Q_{12}$}
\put(153,64){\footnotesize $Q_{22}$}
\put(228,64){\footnotesize $Q_{32}$}

\put(78,114){\footnotesize $Q_{13}$}
\put(153,114){\footnotesize $Q_{23}$}
\put(228,114){\footnotesize $Q_{33}$}

\end{picture}
      \caption{Case $3$ and $5$: 
$\calQ\cap H_8=\{Q_{11}, Q_{22}, Q_{33}\}$.}
\label{fig:case3}
\end{figure}
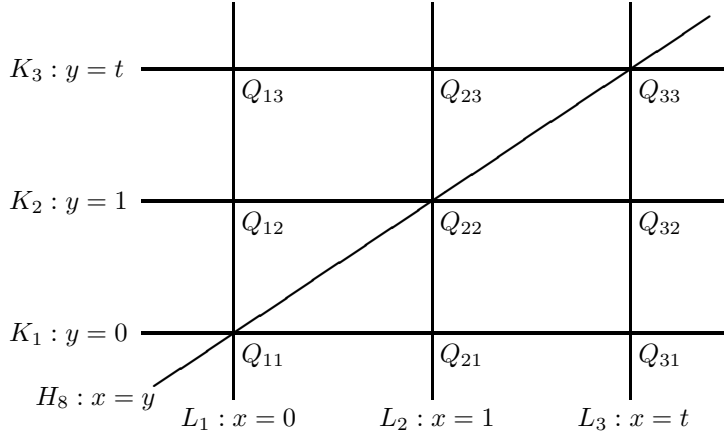

\noindent
{\bf Case 3}: $|H_7\cap \calQ|=2$ and $|H_8\cap \calQ|=3$ 
(Fig. \ref{fig:case3}). 
By changing the numbering of $K_i,L_j$, 
we may assume $H_8\cap \calQ=\{Q_{11},Q_{22},Q_{33}\}$. 
Set $H_7\cap \calQ=\{Q_{i_1j_1}, Q_{i_2j_2}\}$.  
It can be noted that $i_1\neq i_2$ and $j_1\neq j_2$. 
As $\{i_1, i_2\}$ and $\{j_1, j_2\}$ are subsets of 
$\{1,2,3\}$, so the intersection is non-empty. 
Let $k\in \{1,2,3\}$ such that 
$k\in \{i_1,i_2\}\cap \{j_1,j_2\}$. 
Then $H_8\cup K_k\cup L_k$ contains all multiple points of 
$\A$ and $H_8\cap L_k\cap K_k\neq \emptyset$.

%
%
%
%
%

\medskip

\noindent
{\bf Case 4}: $|H_7\cap Q|=|H_8\cap Q|=2$. 

We may assume that $H_8\cap \calQ=\{Q_{11},Q_{22}\}$. 
We can check one-by-one, for any $H_8$, it is 
$\CC_{\leq 3}$ of simple type. 


%
%
%
%
%

\medskip

\noindent
{\bf Case 5}: $|H_7\cap \calQ|=|H_8\cap \calQ|=3$ 
(Fig. \ref{fig:case3}). 
We may assume that $H_8\cap \calQ=\{Q_{11},Q_{22}, Q_{33}\}$. 
We set  $H_7\cap \calQ=\{Q_{1j_1}, Q_{2j_2}, Q_{3j_3}\}$. 
Hence there are six possibilities corresponding to the 
permutation $(j_1,j_2,j_3)$ of $(1,2,3)$. We fix affine 
coordinates as in Figure \ref{fig:case3}. 

%
%
%
%
%

\begin{itemize}
\item[(1)] 
If $(j_1,j_2,j_3)=(1,2,3)$, then $H_7=H_8$. 
\item[(2)] 
If $(j_1,j_2,j_3)=(1,3,2)$. (This implies that $t=-1$.) 
$L_2\cup K_2\cup H_8$ covers all multiple points. 
\item[(3)] 
If $(j_1,j_2,j_3)=(2,1,3)$. (This implies $t=\frac{1}{2}$.) 
$L_1\cup K_1\cup H_8$ covers all multiple points. 
\item[(4)] 
If $(j_1,j_2,j_3)=(3,2,1)$. (This implies $t=2$.) 
$L_1\cup K_1\cup H_8$ covers all multiple points. 
\item[(5)] 
If $(j_1,j_2,j_3)=(3,1,2)$. 
Then $Q_{13}(0,t), Q_{21}(1,0), Q_{32}(1,t)$ are 
collinear if and only if $t=\frac{1\pm\sqrt{-3}}{2}$. 
Hence $\A=\M^\pm$. 
\item[(6)] If $(j_1,j_2,j_3)=(2,3,1)$. Similarly, 
$\A=\M^\pm$. 
\end{itemize}

\endproof

\section{Examples of $9$ and $10$ lines}
\label{sec:ex}

In this section, we will see several examples 
of $9$ and $10$ lines on $\PP^2$ which are 
not covered by previous results. 

\begin{ex}
\label{ex:mc+}
Let $\M^\pm$ be the MacLane arrangement with defining 
equations as in Example \ref{ex:mc}. Consider 
$$
\widetilde{\M}^\pm:=\M^\pm\cup \{H_9\}, 
$$
where $H_9=\{z=0\}$ is the line at infinity 
(Fig. \ref{fig:maclane+}). 

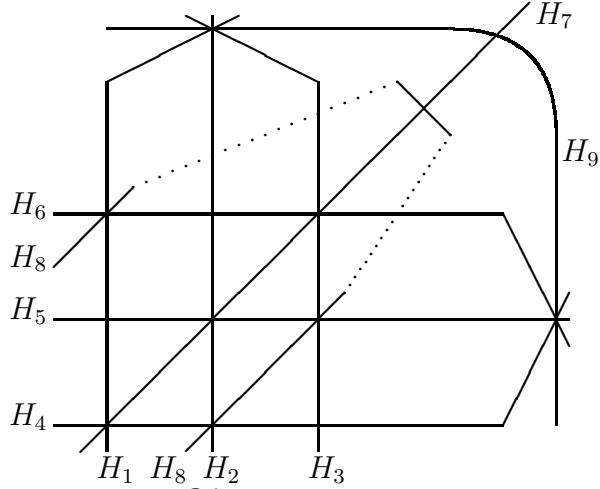
\begin{figure}[htbp]
\begin{picture}(200,170)(0,0)
\thicklines

\put(50,0){\line(0,1){140}}
\put(90,0){\line(0,1){165}}
\put(130,0){\line(0,1){140}}

\put(50,140){\line(2,1){50}}
\put(130,140){\line(-2,1){50}}

\put(30,10){\line(1,0){170}}
\put(30,50){\line(1,0){195}}
\put(30,90){\line(1,0){170}}

\put(200,10){\line(1,2){25}}
\put(200,90){\line(1,-2){25}}

\put(40,0){\line(1,1){170}}

\put(80,0){\line(1,1){60}}
\put(30,70){\line(1,1){30}}
\multiput(60,100)(5,2){21}{\circle*{1}}
\put(160,140){\line(1,-1){20}}
\multiput(140,60)(2,3){21}{\circle*{1}}

\put(46,-10){$H_1$}
\put(86,-10){$H_2$}
\put(126,-10){$H_3$}

\put(13,10){$H_4$}
\put(13,50){$H_5$}
\put(13,90){$H_6$}

\put(212,162){$H_7$}

\put(13,70){$H_8$}
\put(66,-10){$H_8$}

\qbezier(180,160)(220,160)(220,120)
\put(180,160){\line(-1,0){130}}
\put(220,120){\line(0,-1){110}}
\put(222,110){$H_9$}

\end{picture}
      \caption{$\widetilde{\M}^\pm:=\M^\pm\cup \{H_9\}.$}
\label{fig:maclane+}
\end{figure}
The arrangement $\widetilde{\M}^\pm$ is of class $\CC_3$. 
Indeed, all multiple points are contained in 
$H_7\cup H_8\cup H_9$. However since the realization space 
is not connected (Example \ref{ex:mc}), 
it is not $\CC_3$ of simple type. 
\end{ex}

\begin{ex}
\label{ex:fs}
(Falk-Sturmfels arrangements $\FS^\pm$.) 
Let $\gamma_\pm=\frac{1\pm\sqrt{5}}{2}$, and 
define 
$$
\FS^\pm=\{
L_i^\pm, K_i^\pm, H_9^\pm, i=1, 2, 3, 4\}
$$
of $9$ lines as follows (Fig. \ref{fig:fs}): 
$$
\begin{array}{llll}
L_1^\pm: x=0, & 
L_2^\pm: x=\gamma_\pm(y-1), & 
L_3^\pm: y=z, & 
L_4^\pm: x+y=z, \\ 
K_1^\pm: x=z,&
K_2^\pm: x=\gamma_\pm y,&
K_3^\pm: y=0,&
K_4^\pm: x+y=(\gamma_\pm +1)z,\\
H_9^\pm: z=0.&&&
\end{array}
$$

\begin{figure}[htbp]
\begin{picture}(400,214)(20,0)
\thicklines
\put(90,10){\footnotesize $L_1^+$}
\put(100,20){\line(0,1){180}}

\put(7,69){\footnotesize $L_2^+$}
\qbezier(20,70.5573)(100,120)(180,169.443)
\qbezier(180,169.443)(190,175.623)(200,161.803)

\put(7,118){\footnotesize $L_3^+$}
\put(20,120){\line(1,0){160}}
\qbezier(180,120)(190,120)(200,100)

\put(27,183){\footnotesize $L_4^+$}
\put(40,180){\line(1,-1){160}}

\put(110,10){\footnotesize $K_1^+$}
\put(120,20){\line(0,1){160}}
\qbezier(120,180)(120,190)(100,200)

\put(7,48){\footnotesize $K_2^+$}
\qbezier(20,50.5573)(100,100)(200,161.803)

\put(7,98){\footnotesize $K_3^+$}
\put(20,100){\line(1,0){180}}

\put(58,183){\footnotesize $K_4^+$}
\put(72.3607,180){\line(1,-1){107.639}}
\qbezier(180,72.3607)(190,62.3607)(200,20)

\put(15,205){\footnotesize $H_9^+$}
\put(200,0){\line(0,1){190}}
\qbezier(200,190)(200,200)(190,200)
\put(20,200){\line(1,0){170}}


\put(290,10){\footnotesize $L_1^-$}
\put(300,20){\line(0,1){180}}

\put(360,10){\footnotesize $L_2^-$}
\qbezier(262.918,180)(300,120)(361.803,20)
\qbezier(262.918,180)(250.557,200)(238.197,200)

\put(207,118){\footnotesize $L_3^-$}
\put(220,120){\line(1,0){160}}
\qbezier(380,120)(390,120)(400,100)

\put(227,183){\footnotesize $L_4^-$}
\put(240,180){\line(1,-1){160}}

\put(310,10){\footnotesize $K_1^-$}
\put(320,20){\line(0,1){160}}
\qbezier(320,180)(320,190)(300,200)

\put(338,10){\footnotesize $K_2^-$}
\qbezier(238.197,200)(300,100)(350,19.0983)

\put(207,98){\footnotesize $K_3^-$}
\put(220,100){\line(1,0){180}}

\put(227,166){\footnotesize $K_4^-$}
\put(240,167.639){\line(1,-1){140}}
\qbezier(380,27.6393)(390,17.6393)(400,20)

\put(215,205){\footnotesize $H_9^-$}
\put(400,0){\line(0,1){190}}
\qbezier(400,190)(400,200)(390,200)
\put(220,200){\line(1,0){170}}

\end{picture}
      \caption{Falk-Sturmfels arrangements $\FS^+$ and $\FS^-$}
\label{fig:fs}
\end{figure}
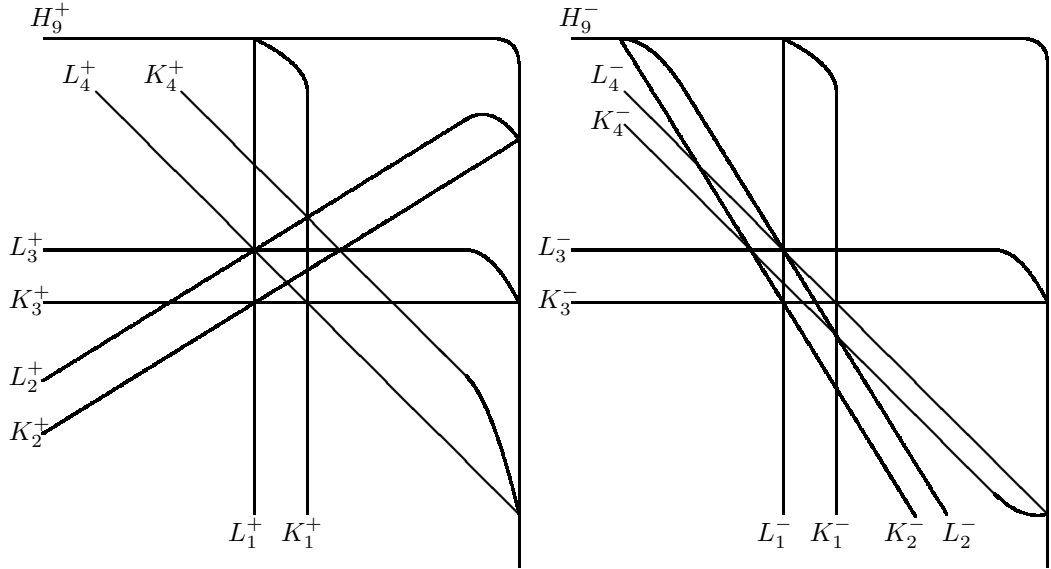
$\FS^+$ and $\FS^-$ have isomorphic incidence relations, which 
are in $\CC_4$ (e.g., multiple points are covered by 
$L_1^\pm\cup L_2^\pm\cup L_3^\pm\cup L_4^\pm$). 
The realization space consists of $2$ connected components 
$\R(I(\FS^\pm))/PGL_3(\C)=\{\FS^+, \FS^-\}$. Thus it is 
the minimal example of $\RR$-realizable arrangement with 
disconnected realization space (Falk-Sturmfels). 
The Galois group action 
$\sqrt{5}\mapsto -\sqrt{5}$ does not induce a 
continuous map of $M(\FS^\pm)$. However there is a 
$PGL_3(\C)$ action 
$(\PP^2, \bigcup_{H\in\FS^+}H)\rightarrow
(\PP^2, \bigcup_{H\in\FS^-}H)$ which maps 
$$
\begin{array}{llll}
L_1^+\longmapsto L_3^-, & 
L_2^+\longmapsto L_4^-, & 
L_3^+\longmapsto L_2^-, & 
L_4^+\longmapsto L_1^-, \\
K_1^+\longmapsto K_3^-, & 
K_2^+\longmapsto K_4^-, & 
K_3^+\longmapsto K_2^-, & 
K_4^+\longmapsto K_1^-, \\
H_9^+\longmapsto H_9^-. 
\end{array}
$$
(In the affine plane the unit square 
$(L_1^+, K_1^+, L_3^+, K_3^+)$ is mapped to the parallelogram 
$(L_3^-, K_3^-, L_2^-, K_2^-)$.) 
In particular, $M(\FS^+)$ and $M(\FS^-)$ are 
homeomorphic and having the isomorphic fundamental 
groups. 
\end{ex}

\begin{ex}
\label{ex:pmi}
(Arrangements $\A^{\pm i}$) 
Define the arrangement 
$$
\A^{\pm i}=\{A_j^\pm, B_j^\pm, C_j^\pm \mid j=1, 2, 3\}, 
$$
of $9$ lines as follows (Fig. \ref{fig:pmi}): 
$$
\begin{array}{lll}
A_1^\pm: x=0, &
A_2^\pm: x=z, &
A_3^\pm: x+y=z, \\
B_1^\pm: y=0, &
B_2^\pm: y=z, &
B_3^\pm: z=0, \\
C_1^\pm: y=\pm\sqrt{-1}x, &
C_2^\pm: y=\mp\sqrt{-1}x+(1\pm\sqrt{-1})z, &
C_3^\pm: x+y=(1\pm\sqrt{-1})z. 
\end{array}
$$
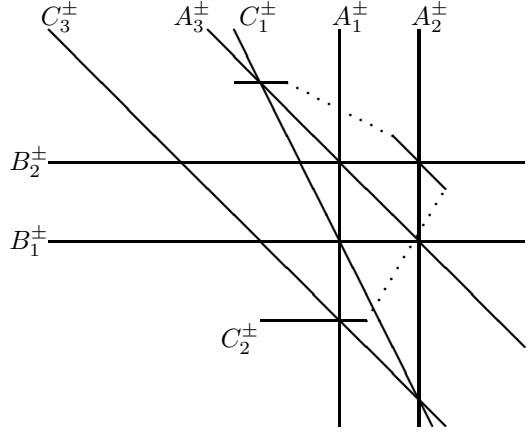
\begin{figure}[htbp]
\begin{picture}(90,160)(0,0)
\thicklines
\put(107, 152){\footnotesize $A_1^\pm$}
\put(110,0){\line(0,1){150}}

\put(137, 152){\footnotesize $A_2^\pm$}
\put(140,0){\line(0,1){150}}

\put(47, 152){\footnotesize $A_3^\pm$}
\put(180,30){\line(-1,1){120}}

\put(-15, 67){\footnotesize $B_1^\pm$}
\put(0,70){\line(1,0){180}}

\put(-15, 97){\footnotesize $B_2^\pm$}
\put(0,100){\line(1,0){180}}


\put(72, 152){\footnotesize $C_1^\pm$}
\put(145,0){\line(-1,2){75}}

\put(65, 30){\footnotesize $C_2^\pm$}
\put(70,130){\line(1,0){20}}
\multiput(90,130)(4,-2){11}{\circle*{1}}
\put(130,110){\line(1,-1){20}}
\multiput(120,40)(2.4,4){13}{\circle*{1}}
\put(80,40){\line(1,0){40}}

\put(-3, 152){\footnotesize $C_3^\pm$}
\put(150,0){\line(-1,1){150}}

\end{picture}
      \caption{$\A^{\pm i}$, where $B_3^\pm$ is the line at infinity.}
\label{fig:pmi}
\end{figure}
It is also in $\CC_4$ (e.g., 
$A_1^\pm\cup A_2^\pm\cup A_3^\pm\cup B_1^\pm$). 
The realization space consists of $2$ connected 
components. As in the case of MacLane arrangement 
(Example \ref{ex:mc}), the complements $M(\A^{\pm i})$ are 
homeomorphic by the complex conjugation. 
\end{ex}

\begin{rk}
\label{rk:9}
Recently the authors verified that, up to $9$ lines, 
these are the complete list of 
disconnected realization spaces. 
Namely, when $|\A|\leq 9$, 
after appropriate re-numbering of $H_1, \dots, H_n$, 
one of the following holds:  
\begin{itemize}
\item[(i)] The realization space 
$\R(I(\A))$ is irreducible (but not necessarily 
$\CC_{\leq 3}$ of simple type, e.g., Pappus 
arrangements), 
\item[(ii)] $\A$ contains the MacLane arrangement $\M^\pm$ 
(Example \ref{ex:mc}, \ref{ex:mc+}), 
\item[(iii)] $\A$ is isomorphic to the Falk-Sturmfels 
arrangement $\FS^\pm$ (Example \ref{ex:fs}), 
\item[(iv)] $\A$ is isomorphic to $\A^{\pm i}$ 
(Example \ref{ex:pmi}). 
\end{itemize}
(Cases (ii), (iii), and (iv) are characterized by 
the minimal field of the realization, 
$\Q(\sqrt{-3}), \Q(\sqrt{5})$, and $\Q(\sqrt{-1})$, 
respectively. It is also concluded from (i) that 
if $I$ is realizable over $\Q$ (with $|\A|\leq 9$), 
then $\R(I)$ is irreducible.) 
The idea of the proof is very similar to 
that of Proposition \ref{prop:8lines} which is 
based on Lemma \ref{lem 3}. 

Consequently, if $I(\A_1)=I(\A_2)$ (with 
$|\A_1|=|\A_2|\leq 9$), $\A_1$ and $\A_2$ are 
transformed to each other by the composition of the 
following operations: 
\begin{itemize}
\item[($\alpha$)] change of numbering, 
\item[($\beta$)] lattice isotopy, 
\item[($\gamma$)] complex conjugation. 
\end{itemize}
In particular, $M(\A_1)$ and $M(\A_2)$ are 
homeomorphic. 
Rybnikov type pairs of arrangements 
require at least $10$ lines. 
\end{rk}

\begin{ex}
\label{ex:extfs}
(Extended Falk-Sturmfels arrangements 
$\widetilde{\FS}^\pm$.) 
Define an arrangement $\widetilde{\FS}^\pm$ of 
$10$ lines by 
adding a line $H_{10}^\pm=\{x=5z\}$ 
to Falk-Sturmfels arrangements 
$\FS^\pm$:  
$$
\widetilde{\FS}^\pm:=\FS^\pm\cup
\{H_{10}^\pm\}. 
$$
$\widetilde{\FS}^\pm$ have the same incidence, 
however there are no ways to transform from 
$\widetilde{\FS}^+$ to 
$\widetilde{\FS}^-$ by operations 
$(\alpha)$, $(\beta)$ and $(\gamma)$. (This fact 
can be proved as follows. First we prove that 
the identity is the only permutation of 
$\{1, \dots, 10\}$ which preserves the incidence. Hence 
if $\widetilde{\FS}^+$ is transformed to $\widetilde{\FS}^-$, 
it sends $L_i^+\mapsto L_i^-, K_i^+\mapsto K_i^-, H_i^+
\mapsto H_i^-$. Deleting $H_{10}^\pm$, 
$\FS^+$ can be transformed to $\FS^-$ with preserving 
the numbering. Note that $\FS^\pm$ are defined over $\RR$ and 
there is no isotopy except for $PGL$ action. There 
should exist a $PGL$ action sending $\FS^+$ to $\FS^-$ which 
preserves the numbering. However it is impossible.) 
The pair $\{\widetilde{\FS}^\pm\}$ is a minimal 
one with such property. 
At this moment the authors do not know 
whether the fundamental groups 
$\pi_1(M(\widetilde{\FS}^\pm))$ are isomorphic. 
\end{ex}

\begin{rk}
We should point out that $\widetilde{\FS}^+$ 
is closer in spirit to examples in \cite[\S $5$]{Ar}. 
\end{rk}

\end{document}